\documentclass[11pt]{article}
\newcommand{\dated}{\mbox{} \hfill {\small [{\tt \today}]}} 
\newtheorem{theorem}{Theorem}[section]
\newtheorem{lemma}[theorem]{Lemma}
\newtheorem{corollary}[theorem]{Corollary}
\newtheorem{proposition}[theorem]{Proposition}
\newtheorem{df}[theorem]{Definition}
\newenvironment{definition}{\begin{df} \rm}{\end{df}}

 \usepackage{amsmath,amssymb,amscd}
\newcommand{\pf}[1]{\trivlist \item[\hskip\labelsep\it #1\ ]}
\newcommand{\varpf}[1]{\trivlist \item[\hskip\labelsep\sc #1:]}
\newcommand{\qedbox}{$\rlap{$\sqcap$}\sqcup$}
\newcommand{\qed}{\qquad \qedbox \endtrivlist}
\newcommand{\varqed}{\hfill \rule{0.6em}{0.6em} \endtrivlist}
\newenvironment{proof}{\pf{Proof}}{\qed}

\newenvironment{remark}{\pf{Remark}}{\endtrivlist}

\newenvironment{examples}{\pf{Examples} 
   \begin{enumerate}}{\end{enumerate} \endtrivlist}
\newenvironment{items}{
  \begin{enumerate} 
                    
  }{\end{enumerate}}
\newenvironment{alphitems}{
  \begin{enumerate} 
                    
  }{\end{enumerate}}
\newenvironment{keywords}{\noindent\small {\it Keywords\/}:}{\vskip 4pt}
\newenvironment{classification}{\noindent\small 2000 {\it Mathematics Subject
Classification\/}:}{\vskip 12pt}

\newcommand{\comps}{{\mathbb C}}

\newcommand{\tensor}{\otimes}

\newcommand{\Tensor}{\hat{\otimes}}

\newcommand{\cstar}{{C^\ast}}

\newcommand{\A}{{\mathfrak A}}
\newcommand{\B}{{\mathfrak B}}

\newcommand{\op}{{\mathrm{op}}}

\newcommand{\VN}{\operatorname{VN}}

\newcommand{\LUC}{\operatorname{\cal LUC}}
\newcommand{\RUC}{\operatorname{\cal RUC}}
\newcommand{\LUCSC}{\operatorname{\cal LUCSC}_0}
\newcommand{\SC}{\operatorname{\cal SC}_0}
\newcommand{\PM}{\operatorname{\it PM}}

\begin{document}
\title{Connes-amenability and normal, virtual diagonals \\ for measure algebras, I}
\author{{\it Volker Runde}\thanks{Research supported by NSERC under grant no.\ 227043-00.}}
\date{}
\maketitle
\begin{center} \large\it
Dedicated to the memory of Barry E.\ Johnson, \\
1937--2002, \\
on whose shoulders many of us stand.
\end{center}
\bigskip\bigskip
\begin{abstract}
We prove that the measure algebra $M(G)$ of a locally compact group $G$ is Connes-amenable if and only if $G$ is amenable.
\end{abstract}
\section*{Introduction}
In \cite{Joh1}, B.\ E.\ Johnson introduced the notion of an amenable Banach algebra, and proved that a locally compact group $G$ is amenable if and
only if its group algebra $L^1(G)$ is amenable. The theory of amenable Banach algebras has been a very active field of research ever since. Once
of the deepest results in this theory is due to A.\ Connes (\cite{Con2}; see also \cite{BP}) and U.\ Haagerup (\cite{Haa}): A $\cstar$-algebra is
amenable if and only if it is nuclear. In \cite{Was1}, S.\ Wassermann showed that a von Neumann algebra is nuclear/amenable if and only if it
is subhomogeneous (see \cite{Run} for a proof that avoids the nuclearity-amenability nexus). This suggests that the definition of amenability from
\cite{Joh1} has to be modified to yield a sufficiently rich theory for von Neumann algebras.
\par
A variant of that definition --- one that takes the dual space structure of a von Neumann algebra into account --- was introduced in 
\cite{JKR}, but is most commonly associated with A.\ Connes' paper \cite{Con1}. For this reason, we refer to this notion of amenability as to
Connes-amenability (the origin of this name seems to be A.\ Ya.\ Helemski\u{\i}'s paper \cite{Hel}). As it turns out, Connes-amenability is the 
``right'' notion of amenability for von Neumann algebras: It is equivalent to several other important properties such as injectivity and semidiscreteness 
(\cite{BP}, \cite{Con1}, \cite{Con2}, \cite{EL}, \cite{Was2}; see \cite[Chapter 6]{LoA} for a self-contained exposition).
\par
The definition of Connes-amenability makes sense for a larger class of Banach algebras (called dual Banach algebras in \cite{Run}). Examples of
dual Banach algebras (other than $W^\ast$-algebras) are: ${\cal B}(E)$, where $E$ is a reflexive Banach space;
$M(G)$, where $G$ is a locally compact group; $\PM_p(G)$, where $p \in (1,\infty)$ and $G$ is a locally compact group (these algebras are called
algebras of $p$-pseudomeasures). The investigation of Connes-amenability for dual Banach algebras which are not $W^\ast$-algebras is still in its
initial stages. Some results on Connes-amenable $W^\ast$-algebras, carry over: For instance, in \cite{LP}, A.\ T.-M.\ Lau and A.\ L.\ T.\ Paterson
showed that, for an inner amenable group $G$, the group von Neumann algebra $\VN(G) =\PM_2(G)$ is Connes-amenable if and only if $G$ is amenable;
this is generalized to $PM_p(G)$ for arbitrary $p \in (1,\infty)$ in \cite{Run}. On the other hand, one cannot expect matters for general dual
Banach algebras to turn out as nicely as for von Neumann algebras: In \cite{Run}, it was shown that ${\cal B}(E)$ is not Connes-amenable if
$E = \ell^p \oplus \ell^q$ with $p,q \in (1,\infty) \setminus \{ 2 \}$ and $p \neq q$.
\par
The dual Banach algebra we are concerned with in this paper is the measure algebra $M(G)$ of a locally compact group $G$. As for von Neumann algebras,
amenability in the sense of \cite{Joh1} is too strong a notion to deal with measure algebras in a satisfactory manner: In \cite{DGH}, H.\ G.\ Dales, F.\ Ghahramani, and
A.\ Ya.\ Helemski\u{\i} prove that $M(G)$ is amenable for a locally compact group $G$ if and only if $G$ is discrete and amenable. In contrast,
Connes-amenability is a much less restrictive demand: Since the amenability of a locally compact group $G$ implies the amenability of $L^1(G)$, and since 
$L^1(G)$ is $w^\ast$-dense in $M(G)$, it follows easily that $M(G)$ is Connes-amenable provided that $G$ is amenable. In this paper, we prove the converse.
\par
I am grateful to S.\ Tabaldyev for discovering a near fatal error in an earlier, stronger version of Lemma \ref{wclem}.
\section{Connes-amenability and normal, virtual diagonals} 
This section is preliminary in character. We collect the necessary definitions we require in the sequel. All of it can 
be found in \cite{Run}, but sometimes our choice of terminology here is different.
\par
Let $\A$ and $\B$ be Banach algebras, and let $E$ be a Banach $\A$-$\B$-bimodule. Then $E^\ast$ becomes a Banach $\B$-$\A$-bimodule via
\[
  \langle x, b \cdot \phi \rangle := \langle x \cdot b, \phi \rangle \quad\text{and}\quad
  \langle x, \phi \cdot a \rangle := \langle a \cdot x, \phi \rangle \qquad (a \in \A, \, b \in \B, \, \phi \in E^\ast, \, x \in E).
\]
\begin{definition}
Let $\A$ and $\B$ be Banach algebras. A Banach $\A$-$\B$-bimodule $E$ is called {\it dual\/} if there is a closed submodule $E_\ast$ of $E^\ast$ such that
$E = (E_\ast)^\ast$.
\end{definition}
\begin{remark}
There is no reason, in general, for $E_\ast$ to be unique. If we refer to the $w^\ast$-topology on a dual Banach module $E$, we always mean
$\sigma(E,E_\ast)$ with respect to a particular, fixed (often obvious) predual $E_\ast$.
\end{remark}
\par
The following definition is due to B.\ E.\ Johnson (\cite{Joh1}):
\begin{definition}
A Banach algebra $\A$ is called {\it amenable\/} if, for every dual Banach $\A$-bimodule $E$, every bounded derivation $D \!: \A \to E$ is inner.
\end{definition}
\par
We are interested in a particular class of Banach algebras:
\begin{definition} \label{ddef}
A Banach algebra $\mathfrak A$ is called {\it dual\/} if it is dual as a Banach $\A$-bimodule.
\end{definition}
\begin{remark}
A Banach algebra which is also a dual space is a dual Banach algebra if and only if multiplication is separately $w^\ast$-continuous.
\end{remark}      
\begin{examples}
\item Every $W^\ast$-algebra is dual.
\item If $E$ is a reflexive Banach space, then ${\cal B}(E) = (E \Tensor E^\ast)^\ast$ is dual.
\item If $G$ is a locally compact group, then $M(G) = {\cal C}_0(G)^\ast$ is dual.
\item If $\A$ is an Arens regular Banach algebra, then $\A^{\ast\ast}$ is dual.
\end{examples}
\par
The following choice of terminology is motivated by the von Neumann algebra case:
\begin{definition}
Let $\A$ and $\B$ be dual Banach algebras. A dual Banach $\A$-$\B$-bimodule is called {\it normal\/} if, for each $x \in E$, the maps
\[
  \A \to E, \quad a \mapsto a \cdot x 
\]
and
\[
  \B \to E, \quad b \mapsto x \cdot b
\]
are $w^\ast$-continuous.
\end{definition}
\par
We can now define Connes-amenable, dual Banach algebras:
\begin{definition}
A dual Banach algebra $\A$ is called {\it Connes-amenable\/} if, for every normal, dual Banach $\A$-bimodule $E$, every bounded, $w^\ast$-continuous
derivation $D \!: \A \to E$ is inner.
\end{definition}
\par
Amenability in the sense of \cite{Joh1}, can be intrinsically characterized in terms of so-called approximate and virtual diagonals (\cite{Joh2}).
There is a related notion for Connes-amenable, dual Banach algebras.
\par
If $E_1, \ldots, E_n$ and $F$ are dual Banach spaces, we write ${\cal L}_{w^\ast}(E_1, \ldots, E_n ; F)$ for the bounded, separately
$w^\ast$-continuous, $n$-linear maps from $E_1 \times \cdots \times E_n$ into $F$. In case $E_1 = \cdots = E_n =: E$, we simply let
${\cal L}_{w^\ast}^n(E,F) := {\cal L}_{w^\ast}(E_1, \ldots, E_n ; F)$.
\par
Let $\A$ and $\B$ be Banach algebras. Then $\A \Tensor \B$ becomes a Banach $\A$-$\B$-bimodule via
\begin{equation} \label{tensormod}
  a\cdot (x \tensor y ) := ax \tensor y \quad\text{and}\quad (x \tensor y) \cdot b := x \tensor y b \qquad (a,x \in \A, \, b,y \in \B).
\end{equation}
Suppose that $\A$ and $\B$ are dual. It is then routinely checked that ${\cal L}_{w^\ast}^2(\A,\B;\comps)$ is a closed $\B$-$\A$-submodule of
$(\A \Tensor \B)^\ast$. 
\par
Let $\A$ be a dual Banach algebra, and let $\Delta_\A \!: \A \Tensor \A \to \A$ denote the diagonal operator induced by $\A \times \A \ni (a,b) \mapsto ab$.
Since multiplication in $\A$ is separately $w^\ast$-continuous, we have $\Delta^\ast_\A \A_\ast \subset {\cal L}_{w^\ast}^2(\A,\comps)$. Taking the adjoint of 
$\Delta^\ast_\A |_{\A_\ast}$, we may thus extend $\Delta_\A$ to ${\cal L}_{w^\ast}^2(\A,\comps)^\ast$ as an $\A$-bimodule homomorphism (we denote this extension by 
$\Delta_{w^\ast}$).
\begin{definition}
Let $\A$ be a dual Banach algebra. An element ${\rm M} \in {\cal L}_{w^\ast}^2(\A,\comps)^\ast$ is called a {\it normal, virtual diagonal\/} for $\A$ if
\[
  a \cdot {\rm M} = {\rm M} \cdot a \quad\text{and}\quad a \Delta_{w^\ast} {\rm M} = a \qquad (a \in \A). 
\]
\end{definition} 
\par
One connection between Connes-amenability and the existence of normal, virtual diagonals is fairly straightforward (\cite{Eff}, \cite{CG}): If
$\A$ has a normal, virtual diagonal, then $\A$ is Connes-amenable; in fact, it implies a somewhat stronger property (\cite{Run}). The main problem
with proving the converse is that, in general, the dual module ${\cal L}_{w^\ast}^2(\A,\comps)^\ast$ need not be normal. For von Neumann algebras, however,
Connes-amenability and the existence of normal, virtual diagonals are even equivalent (\cite{Eff}, \cite{EK}). We suspect, but have been unable to prove --- except in the discrete
case --- that the same is true for the measure algebras of locally compact groups.
\section{Separately ${\cal C}_0$-functions on locally compact Hausdorff spaces}
Our notation is standard: For a topological space $X$, we write ${\cal C}_b(X)$ for the bounded, continuous functions on $X$; if $X$ is locally compact and
Hausdorff, ${\cal C}_0(X)$ (or rather ${\cal C}(X)$ if $X$ is compact) denotes the continuous functions on $X$ vanishing at infinity, and $M(X) \cong {\cal C}_0(X)^\ast$ 
is the space of regular Borel measures on $X$.
\par
Let $X$ and $Y$ be locally compact Hausdorff spaces. In this section, we give a description of ${\cal L}_{w^\ast}(M(X),M(Y);\comps)$ as a space of separately
continuous functions on $X \times Y$.
\begin{definition} \label{csdef}
Let $X$ and $Y$ be locally compact Hausdorff spaces. A bounded function $f \!: X \times Y \to \comps$ is called {\it separately ${\cal C}_0$\/} if:
\begin{alphitems}
\item for each $x \in X$, the function
\[
  Y \to \comps, \quad y \mapsto f(x,y) 
\]
belongs to ${\cal C}_0(Y)$;
\item for each $y \in Y$, the function
\[
  X \to \comps, \quad x \to f(x,y)
\]
belongs to ${\cal C}_0(X)$.
\end{alphitems}
We define $\SC(X \times Y)$ as the collection of all separately ${\cal C}_0$-functions.
\end{definition}
\begin{lemma} \label{sepc0l1}
Let $X$ and $Y$ be locally compact Hausdorff spaces, and let $f \in \SC(X \times Y)$. Then the following hold:
\begin{items}
\item for each $\mu \in M(X)$, the function
\[
  Y \to \comps, \quad y \mapsto \int_X f(x,y) \, d\mu(x)
\]
belongs to ${\cal C}_0(Y)$;
\item for each $\nu \in M(Y)$, the function
\[
  X \to \comps, \quad x \mapsto \int_Y f(x,y) \, d\nu(y)
\]
belongs to ${\cal C}_0(X)$.
\end{items}
\end{lemma}
\begin{proof}
We only prove (i).
\par
Let $\mu \in M(X)$. Since the measures with compact support are norm dense in $M(X)$, there is no loss of generality if we suppose that $X$ is
compact. Suppose that $Y$ is not compact (the compact case is easier), and let $Y_\infty$ be its one-point-compactification. Extend $f$
to $X \times Y_\infty$ by letting
\[
  f(x,\infty) = 0 \qquad (x \in X),
\]
so that $f$ is separately continuous on $X \times Y_\infty$. Let $\tau$ be the topology of pointwise convergence on ${\cal C}(X)$. Since
the map
\[
  Y_\infty \to {\cal C}(X), \quad y \mapsto f(\cdot,y)
\]
is continuous with respect to the given topology on $Y_\infty$ and to $\tau$ on ${\cal C}(X)$, the set
\[
  K := \{ f(\cdot,y) : y \in Y_\infty \} 
\]
is $\tau$-compact. By \cite[Th\'eor\`eme 5]{Gro}, this means that $K$ is weakly compact, so that the weak topology and $\tau$ coincide on $K$. 
Let $(y_\alpha)_\alpha$ be a convergent net in $Y_\infty$ with limit $y$. Since $f(\cdot,y_\alpha) \stackrel{\tau}{\to} f(\cdot,y)$,
it follows that
\[
  \lim_\alpha \int_X f(x,y_\alpha) \, d\mu(x) = \int_X f(x,y) \,d \mu(x).
\]
This means that the function
\[
  Y_\infty \to \comps, \quad y \mapsto \int_X f(x,y) \,d \mu(x)
\]
is continuous on $Y_\infty$; since it vanishes at $\infty$ by definition, this establishes (i).
\end{proof}
\begin{remark}
For compact spaces, Lemma \ref{sepc0l1} is well known (\cite[Theorem A.20]{Bur}).
\end{remark}
\begin{lemma} \label{sepc0l2}
Let $X$ and $Y$ be locally compact Hausdorff spaces, and let $f \in \SC(X \times Y)$. Then the bilinear map
\begin{equation} \label{eq1}
  \Phi_f \!: M(X) \times M(Y) \to \comps, \quad (\mu,\nu) \mapsto \int_Y \int_X f(x,y) \, d\mu(x) \, d\nu(y)
\end{equation}
belongs to ${\cal L}_{w^\ast}(M(X),M(Y); \comps)$.
\end{lemma}
\begin{proof}
Clearly, $\Phi_f$ is bounded, and it is immediate from Lemma \ref{sepc0l1}(i) that it is $w^\ast$-continuous in the second variable.
Since $f$ is separately continuous, it is $\mu \tensor \nu$-measurable for all $\mu \in M(X)$ and $\nu \in M(Y)$ by \cite{Joh0}. It follows that the integral
in (\ref{eq1}) not only exists, but --- by Fubini's theorem --- is independent of the order of integration, i.e.\
\[
  \Phi_f(\mu,\nu) = \int_X \int_Y f(x,y) \, d\nu(y) \, d\mu(x) \qquad (\mu \in M(X), \, \nu \in M(Y)).
\]
It then follows from Lemma \ref{sepc0l1}(ii) that $\Phi_f$ is also $w^\ast$-continuous in the first variable.
\end{proof}
\begin{proposition} \label{sepc0prop}
Let $X$ and $Y$ be locally compact Hausdorff spaces. Then
\begin{equation} \label{eq2}
  \SC(X \times Y) \to {\cal L}^2_{w^\ast}(M(X), M(Y); \comps), \quad f \mapsto \Phi_f
\end{equation}
is an isometric isomorphism.
\end{proposition}
\begin{proof}
Clearly, $\| \Phi_f \| \leq \| f \|_\infty$ for all $f \in \SC(X \times Y)$. On the other hand,
\begin{eqnarray*}
  \| \Phi_f \| & \geq & \sup \{ | \Phi_f(\delta_x,\delta_y)| : x\in X, \, y \in Y \} \\
  & = & \sup \{ | f(x,y)| : x\in X, \, y \in Y \} \\
  & = & \| f \|_\infty \qquad (f \in \SC(X \times Y)),
\end{eqnarray*}
so that (\ref{eq2}) is an isometry.
\par
Let $\Phi \in {\cal L}_{w^\ast}(M(X),M(Y);\comps)$ be arbitrary, and define
\[
  f \!: X \times Y \to \comps, \quad (x,y) \mapsto \Phi(\delta_x, \delta_y).
\]
It is immediate that $f \in \SC(X \times Y)$ such that $\Phi_f(\delta_x,\delta_y) = \Phi(\delta_x,\delta_y)$ for all $x \in X$ and
$y \in Y$. Separate $w^\ast$-continuity yields that $\Phi = \Phi_f$.
\end{proof}
\par
We shall, from now on, identify $\SC(X \times Y)$ and ${\cal L}_{w^\ast}(M(X),M(Y);\comps)$ as Banach spaces.
\begin{proposition} \label{sepc0prop2}
Let $X$ and $Y$ be locally compact Hausdorff spaces. Then the map 
\[
  M(X \times Y) \to {\cal L}^2_{w^\ast}(M(X), M(Y); \comps)^\ast, \quad \mu \mapsto \Psi_\mu, 
\]
where
\begin{equation} \label{inteq}
  \Psi_\mu(f) := \int_{X \times Y} f(x,y) \, d\mu(x,y) \qquad (\mu \in M(X \times Y), \, f \in  \SC(X \times Y)),
\end{equation} 
is an isometry.
\end{proposition}
\begin{proof}
By \cite{Joh0} again, the integral in (\ref{inteq}) is well-defined. Since ${\cal C}_0(X \times Y) \subset \SC(X \times Y)$,
it follows at once that $\| \Psi_\mu \| = \| \mu \|$ holds for all $\mu \in M(X \times Y)$.
\end{proof}
\section{Separately ${\cal C}_0$-functions on locally compact groups}
Let $G$ and $H$ be locally compact groups. Then $\SC(G \times H)$ becomes a Banach $M(H)$-$M(G)$-bimodule through the following
convolution formulae for $f \in \SC(G \times H)$, $\mu \in M(H)$, and $\nu \in M(G)$:
\[
  (\mu \cdot f)(g,h) := \int_H f(g, hk) \, d\mu(k) \qquad (g \in G, \, h \in H)
\]
and
\[
  (f \cdot \nu)(g,h) := \int_G f(kg, h) \, d\nu(k) \qquad (g \in G, \, h \in H).
\]
\par
The following extension of Proposition \ref{sepc0prop} is then routinely checked:
\begin{proposition} \label{modprop}
Let $G$ and $H$ be locally compact groups. Then
\begin{equation} \label{modhom0}
  \SC(G \times H) \to {\cal L}^2_{w^\ast}(M(G), M(H); \comps), \quad f \mapsto \Phi_f
\end{equation}
as defined in Proposition\/ {\rm \ref{sepc0prop}} is an isometric isomorphism of Banach $M(H)$-$M(G)$-bi\-mod\-ules.
\end{proposition}
\begin{proof}
The $M(H)$-$M(G)$-bimodule action on $\SC(G \times H)$ induces an $M(G)$-$M(H)$-bi\-mod\-ule action on $\SC(G \times H)^\ast$. Embedding
$M(G) \Tensor M(H)$ into $\SC(G \times H)^\ast$, we need to show that $M(G) \Tensor M(H)$ is a $M(G) \Tensor M(H)$-submodule of $\SC(G \times H)^\ast$
such that the module actions are the canonical ones (see (\ref{tensormod})).
\par
Let $\kappa, \mu \in M(G)$, and let $\nu \in M(H)$. Then we have for $f \in \SC(G \times H)$:
\begin{eqnarray*}
  \langle f, \kappa \cdot (\mu \tensor \nu) \rangle & = & \langle f \cdot \kappa, \mu \tensor \nu \rangle \\
  & = & \int_H \int_G \int_G f(kg,h) \, d\kappa(k) \, d\mu(g) \, d\nu(h) \\
  & = & \int_H \int_G f(g,h) \, d(\kappa \ast \mu)(g) \, d\nu(h) \\
  & = & \langle f, \kappa \ast \mu \tensor \nu \rangle.
\end{eqnarray*}
\par
An analogous property holds for the right $M(H)$-module action on $\SC(G \times H)^\ast$.
\end{proof}
\begin{remark}
It is easy to see that ${\cal C}_0(G \times H)$ is a closed $M(H)$-$M(G)$-submodule of $\SC(G \times H)$, so that $M(G \times H) \cong {\cal C}_0(G \times H)^\ast$
is a quotient of $\SC(G \times H)^\ast$. It is easily checked that
\[
  \mu \cdot \nu = (\mu \tensor \delta_e) \ast \nu \qquad (\mu \in M(G), \, \nu \in M(G \times H))
\]
and
\[
  \nu \cdot \mu = \nu \ast (\delta_e \tensor \mu) \qquad (\mu \in M(H), \, \nu \in M(G \times H)).
\]
\end{remark}
\par
We have the following:
\begin{proposition} \label{modprop2}
Let $G$ and $H$ be locally compact groups. Then:
\begin{items}
\item $M(G \times H)$ is a normal, dual Banach $M(G)$-$M(H)$-bimodule.
\item The map 
\begin{equation} \label{modhom}
  M(G \times H) \to {\cal L}^2_{w^\ast}(M(G), M(H); \comps)^\ast, \quad \mu \mapsto \Psi_\mu,
\end{equation}
as defined in Proposition\/ {\rm \ref{sepc0prop2}}, is an isometric homomorphism of Banach $M(G)$-$M(H)$-bimodules.
\end{items}
\end{proposition}
\begin{proof}
The maps
\[
  M(G) \to M(G \times H), \quad \mu \mapsto \mu \tensor \delta_e
\]
and
\[
  M(H) \to M(G \times H), \quad \nu \mapsto \delta_e \tensor \nu
\]
are $w^\ast$-continuous. In view of the preceding remark and the fact that $M(G \times H)$ is a dual Banach algebra, (i) is immediate.
\par
For (ii), let $\mu \in M(G)$ and $\nu \in M(G \times H)$. Then we have for any $f \in \SC(G \times H)$:
\begin{eqnarray*}
   \langle f, \mu \cdot \Psi_\nu \rangle & = & \langle f \cdot \mu, \Psi_\nu \rangle \\
   & = & \int_{G \times H} \int_G f(kg,h) \, d\mu(k) \, d\nu(g,h) \\
   & = & \int_{G \times H} \int_{G \times H} f(kg,k'h) \, d(\mu \tensor \delta_e)(k,k') \, d\nu(g,h) \\
   & = & \int_{G \times H} f(g,h) \, d((\mu \tensor \delta_e) \ast \nu)(g,h) \\
   & = & \int_{G \times H} f(g,h) \, d(\mu \cdot \nu)(g,h) \\ 
   & = & \langle f, \Psi_{\mu \cdot \nu} \rangle.
\end{eqnarray*}
Hence, (\ref{modhom}) is a left $M(G)$-module homomorphism. 
\par
Analogously, one shows that (\ref{modhom}) is a right $M(H)$-module homomorphism.
\end{proof}
\par
With these preparations made, we can already give an alternative proof of \cite[Proposition 5.2]{Run}. 
\par
For any locally compact group $G$, the operator $\Delta_\ast := \Delta^\ast_{M(G)} |_{{\cal C}_0(G)}$ is given by
\[
  (\Delta_\ast f)(g,h) = f(gh) \qquad (f \in {\cal C}_0(G), \, g,h \in G).
\]
If $G$ is compact, $\Delta_\ast$ maps ${\cal C}_0(G) = {\cal C}(G)$ into ${\cal C}(G \times G) = {\cal C}_0(G \times G)$. Hence, $\Delta_{w^\ast}$
drops to an $M(G)$-bimodule homomorphism $\Delta_{0,w^\ast} \!: M(G \times G) \to M(G)$.
\begin{proposition} \label{compact}
Let $G$ be a compact group. Then there is a normal, virtual diagonal for $M(G)$.
\end{proposition}
\begin{proof}
Since $G$ is amenable, $M(G)$ is Connes-amenable (this is the easy direction of Theorem \ref{thm}).
\par
Define a $w^\ast$-continuous derivation
\[
  D \!: M(G) \to M(G \times G), \quad \mu \mapsto \mu \tensor \delta_e - \delta_e \tensor \mu.
\]
It is immediate that $D$ attains its values in $\ker \Delta_{0,w^\ast}$. Being the kernel of a $w^\ast$-con\-ti\-nu\-ous $M(G)$-bimodule homomorphism between
normal, dual Banach $M(G)$-bimodules, $\ker \Delta_{0,w^\ast}$ is a normal, dual Banach $M(G)$-bimodule in its own right. Since $M(G)$ is Connes-amenable, 
there is ${\rm N} \in \ker \Delta_{0,w^\ast}$ such that
\[
  D\mu = \mu \cdot {\rm N} - {\rm N} \cdot \mu \qquad (\mu \in M(G)).
\]
Letting ${\rm M} := \delta_e \tensor \delta_e - {\rm N}$, and embedding ${\rm M}$ into $\SC(G \times G)^\ast$ via Proposition \ref{modprop2}, we obtain a 
normal, virtual diagonal for $M(G)$.
\end{proof}
\begin{remark}
The proof of Proposition \ref{compact}, does not carry over to non-compact, locally compact groups with Connes-amenable measure algebra because,
for non-compact $G$, we no longer have $\Delta_\ast {\cal C}_0(G) \subset {\cal C}_0(G \times G)$; in fact, it is easy to see that $\Delta_\ast{\cal C}_0(G)
\cap {\cal C}_0(G \times G) = \{ 0 \}$ whenever $G$ is not compact.
\end{remark}
\section{A left introverted subspace of separately ${\cal C}_0$-functions}
For general, possibly non-compact, locally compact groups, we need a Banach $M(G)$-bimodule that can play the r\^ole of $M(G \times G)$ in the proof
of Proposition \ref{compact}.
\par
Let $G$ be a locally compact group. For a function $f \!: G \to \comps$ and for $g \in G$, define functions $L_g f , R_g f \!: G \to \comps$ through
\[
  (L_g f)(h) := f(gh) \quad\text{and}\quad (R_g f)(h) := f(hg) \qquad (h \in G).
\]
A closed subspace $E$ of $\ell^\infty(G)$ is called {\it left invariant\/} if $L_g f \in E$ for each $f \in E$ and $g \in G$. A left invariant subspace $E$
of $\ell^\infty(G)$ is called {\it left introverted\/} if, for each $\phi \in E^\ast$, the function
\[
  \phi \bullet f \!: G \to \comps, \quad g \mapsto \langle L_g f , \phi \rangle 
\]
belongs again to $E$.
\begin{examples}
\item $\ell^\infty(G)$ is trivially left introverted.
\item ${\cal C}_0(G)$ is left introverted (\cite[(19.5) Lemma]{HR}). 
\item The space
\[
  \LUC(G) := \{ f \in {\cal C}_b(G) : \text{$G \ni g \mapsto L_g f$ is norm continuous} \}
\]
of {\it left uniformly continuous functions\/} on $G$ is left introverted (\cite[(2.11) Proposition]{Pat}). 
\end{examples}
\par
If $E$ is a left introverted subspace of $\ell^\infty(G)$, then $E^\ast$ is a Banach algebra in a natural manner:
\[
  \langle \phi \ast \psi,  f \rangle := \langle \psi \bullet f, \phi \rangle \qquad (\phi, \psi \in E^\ast, \, f \in E).
\]
In case $E = {\cal C}_0(G)$, this is the usual convolution product on $M(G)$.
\par
We now define a certain space of separately ${\cal C}_0$-functions which is, as we shall see, left introverted. For any locally compact group $G$, let $G_{\LUC}$ denote the character space of the commutative $\cstar$-algebra 
$\LUC(G)$. The multiplication $\ast$ on $\LUC(G)^\ast$ turns $G_{\LUC}$ into a compact semigroup with continuous right multiplication that contains $G$ as a dense subsemigroup
(\cite{BJM}). Also, we use $G^\op$ to denote the same group, but with reversed multiplication.
\begin{definition} \label{lucsdef}
For locally compact groups $G$ and $H$, let
\begin{eqnarray*}
  \lefteqn{\LUCSC(G \times H^\op)} & & \\
  & := & \{ f \in \LUC(G \times H^\op) : \text{$\omega \bullet f \in \SC(G \times H)$ for all $\omega \in (G \times H^\op)_{\LUC}$} \}.
\end{eqnarray*}
\end{definition}
\begin{remark}
If both $G$ and $H$ are compact, then $\LUCSC(G \times H^\op) = {\cal C}(G \times H)$.
\end{remark}
\begin{lemma} \label{wclem}
Let $G$ and $H$ be locally compact groups, let $f \in \LUCSC(G \times H^\op)$, and let $h \in H$. Then $\{ L_{(g,h)} f : g \in G \}$
is relatively weakly compact. 
\end{lemma}
\begin{proof}
The claim is clear for compact $G$, so that we may suppose without loss of generality that $G$ is not compact.
\par
By \cite[Th\'eor\`eme 5]{Gro}, it is sufficient to show that $\{ L_{(g,h)} f : g \in G \}$ is relatively compact in $\LUC(G \times H^\op)$ with respect to the 
topology of pointwise convergence on $(G \times H^\op)_{\LUC}$. Also, we may suppose without loss of generality that $h = e$.
\par
Let $\hat{f} \in {\cal C}((G \times H^\op)_{\LUC})$ denote the Gelfand transform of $f$. The map
\[
  G \to \comps, \quad g \mapsto \hat{f}((\delta_g \tensor \delta_e) \ast \omega)
\]
is continuous for each $\omega \in (G \times H^\op)_{\LUC}$. Let $G_\infty$ denote the one-point-compactification of $G$. Let $( g_\alpha )_\alpha$ be a net in $G$ with $g_\alpha \to \infty$. For any $\omega 
\in (G \times H^\op)_{\LUC}$, we then have
\[
  \hat{f}((\delta_{g_\alpha} \tensor \delta_e) \ast \omega) =(\omega \bullet f)(g_\alpha, e) \to 0
\]
because $\omega \bullet f \in \SC(G \times H)$. Hence,
\begin{equation} \label{lefttransl}
  G \to \LUC(G \times H^\op), \quad g \mapsto L_{(g,e)}f
\end{equation}
extends as a continuous map to $G_\infty$, where $\LUC(G \times H^\op)$ is equipped with the topology of pointwise convergence on $(G \times H^\op)_{\LUC}$. As the continuous image of
the compact space $G_\infty$, the range of (\ref{lefttransl}) is compact in the topology of pointwise convergence on $(G \times H^\op)_{\LUC}$.
\end{proof}
\begin{proposition} \label{lucprop0}
Let $G$ and $H$ be locally compact groups. Then $\LUCSC(G \times H^\op)$ is left introverted.
\end{proposition}
\begin{proof}
Let $f \in \LUCSC(G \times H^\op)$, and let $\phi \in  \LUCSC(G \times H^\op)^\ast$. Since $\LUC(G \times H^\op)$ is left introverted, it is immediate that
$\phi \bullet f \in \LUC(G \times H^\op)$. 
\par
We first claim that $\phi \bullet f \in \SC(G \times H)$.
\par
Fix $h \in H$; we will show that $(\phi \bullet f)( \cdot, h)$, i.e.\ the function
\[
  G \to \comps,\quad g \mapsto \langle L_{(g,h)} f , \phi \rangle 
\]
belongs to ${\cal C}_0(G)$. Since $(\phi \bullet f)( \cdot, h)$ is clearly continuous, all we have to show is that it vanishes at $\infty$. Suppose without loss of generality that
$G$ is not compact, and let $(g_\alpha )_\alpha$ be a net in $G$ such that $g_\alpha \to \infty$. Let $\tau$ denote the topology of pointwise convergence on $G \times H$. 
It is clear that $L_{(g_\alpha, h)} f \stackrel{\tau}{\to} 0$. Since $\{ L_{(g,h)} f : g \in G \}$ is relatively weakly compact by Lemma \ref{wclem}, the weak 
topology and $\tau$ coincide on the weak closure of $\{ L_{(g,h)} f : g \in G \}$, so that, in particular, $\langle L_{(g_\alpha, h)} f , \phi \rangle \to 0$.
\par
Analogously, one sees that  $(\phi \bullet f)(g, \cdot) \in {\cal C}_0(H)$ for each $g \in G$.
\par
Let $\omega \in (G \times H^\op)_{\LUC}$. Since, by the foregoing,
\[
  \omega \bullet (\phi \bullet f) = (\omega \ast \phi) \bullet f \in \SC(G \times H),
\]
it follows that $\phi \bullet f \in \LUCSC(G \times H^\op)$.
\end{proof}
\begin{theorem} \label{ucthm}
Let $G$ and $H$ be locally compact groups. Then we have:
\begin{items}
\item $\LUCSC(G \times H^\op)$ is a closed $M(H)$-$M(G)$-submodule of $\SC(G \times H^\op)$.
\item $\LUCSC(G \times H^\op)^\ast$ is a normal, dual Banach $M(G)$-$M(H)$-bimodule.
\item If $H = G$, then $\Delta_\ast$ maps ${\cal C}_0(G)$ into $\LUCSC(G \times G^\op)$.
\end{items}
\end{theorem}
\begin{proof}
For (i), first note that it is routinely checked that $\mu \cdot f, f \cdot \nu \in \LUC(G \times H^\op)$ for all $f \in \LUCSC(G \times H^\op)$ and all $\mu \in M(G)$ and $\nu \in M(H)$. Fix $f \in 
\LUCSC(G \times H^\op)$, $\mu \in M(G)$, $\nu \in M(H)$, and let $\omega \in (G \times H^\op)_{\LUC}$. Since
\[
  \omega \bullet (\mu \cdot f)(g,h) = \langle \mu \cdot L_{(g,h)} f, \omega \rangle \qquad ((g,h) \in G \times H^\op),
\]
an application of Lemma \ref{wclem} as in the proof of Proposition \ref{lucprop0} yields that $\omega \bullet (\mu \cdot f) \in \SC(G \times H^\op)$. A similar, but easier argument yields that
$\omega \bullet (f \cdot \nu) \in \SC(G \times H^\op)$.
\par
For (ii), first observe that the canonical embedding of $M(G \times H^\op)$ into $\LUC(G \times H^\op)^\ast$ via integration is an algebra homomorphism.
If we view $M(G \times H^\op)$ as a $M(G)$-$M(H)$-submodule of $\LUCSC(G \times H^\op)^\ast$ (through Proposition \ref{modprop2}(ii)), we see routinely that
\begin{equation} \label{conv5}
  \mu \cdot \nu = (\mu \tensor \delta_e) \ast \nu |_{\LUCSC(G \times H)} \qquad (\mu \in M(G), \, \nu \in M(G \times H^\op)).
\end{equation}
Fix $\mu \in M(G)$. By (the simple direction of) \cite{Lau} --- actually already proven in \cite{Won} ---, the map
\begin{equation} \label{tony}
  \LUC(G \times H^\op)^\ast \to \LUC(G \times H^\op)^\ast, \quad \phi \mapsto (\mu \tensor \delta_e) \ast \phi
\end{equation}
is $w^\ast$-continuous. Let $\phi \in \LUC(G \times H^\op)^\ast$ be arbitrary, and choose a net $( \nu_\alpha )_\alpha$ in $M(G \times H^\op)$ that
converges to $\phi$ in the $w^\ast$-topology (the existence of such a net follows with a simple Hahn--Banach argument). Then (\ref{conv5}) and the 
$w^\ast$-continuity of (\ref{tony}), yield that
\[
  \mu \cdot \phi = \text{$w^\ast$-}\lim_\alpha \mu \cdot \nu_\alpha = \text{$w^\ast$-}\lim_\alpha (\mu \tensor \delta_e) \ast \nu_\alpha
  = (\mu \tensor \delta_e) \ast \phi.
\]
Let $f \in \LUCSC(G \times H^\op)$. Then we have
\begin{equation} \label{modop}
  \langle f, \mu \cdot \phi \rangle = \langle f, (\mu \tensor \delta_e) \ast \phi \rangle = \langle \phi \bullet f , \mu \tensor \delta_e \rangle.
\end{equation}
By Proposition \ref{lucprop0}, $\LUCSC(G \times H^\op)$ is left introverted, so that, in particular, $\phi \bullet f \in \SC(G \times H^\op)$. Let $(\mu_\alpha )_\alpha$ be a net in $M(G)$ that converges to $\mu$ in the 
$w^\ast$-topology. Then (\ref{modop}) yields:
\begin{eqnarray*}
  \lim_\alpha \langle f, \mu_\alpha \cdot \phi \rangle & = & \lim_\alpha \langle \mu_\alpha \tensor \delta_e, \phi \bullet f \rangle \\
  & = & \lim_\alpha \langle (\phi \bullet f )(\cdot, e), \mu_\alpha \rangle \\
  & = & \langle (\phi \bullet f )(\cdot, e), \mu \rangle \\
  & = & \langle \phi \bullet f,  \mu \tensor \delta_e \rangle \\
  & = & \langle f, \mu \cdot \phi \rangle.
\end{eqnarray*}
\par
It follows that, for any $\phi \in \LUCSC(G \times H^\op)^\ast$, the map
\[
  M(G) \to \LUCSC(G \times H^\op)^\ast, \quad \mu \mapsto \mu \cdot \phi
\]
is $w^\ast$-continuous. Noting that
\[
  \phi \cdot \nu = \phi \ast (\delta_e \tensor \nu) \qquad (\nu \in M(H)),
\]
we see analoguously that
\[
  M(H) \to \LUCSC(G \times H^\op)^\ast, \quad \nu \mapsto \phi \cdot \nu
\]
is $w^\ast$-continuous for all $\phi \in \LUCSC(G \times H^\op)^\ast$. This proves (ii).
\par
Suppose that $H = G$. It is well known ${\cal C}_0(G) \subset \LUC(G) \cap \RUC(G)$, where
\[
  \RUC(G) := \{ f \in {\cal C}_b(G) : \text{$G \ni g \mapsto R_g f$ is norm continuous} \}.
\]
Let $f \in {\cal C}_0(G)$, and note that
\[
  L_{(g,h)} \Delta_\ast f = \Delta_\ast(L_g R_h f) \qquad ((g,h) \in G \times G^\op).
\]
The norm continuity of $\Delta_\ast$ shows that $\Delta_\ast f \in \LUC(G \times G^\op)$. To show that $\Delta_\ast f \in \LUCSC(G \times G^\op)$, let $\omega \in (G \times G^\op)_{\LUC}$. Let $((g_\alpha, h_\alpha))_\alpha$ be a
net in $G \times G^\op$ such that $(g_\alpha, h_\alpha) \to \omega$. Passing to a subnet, we may suppose that $( g_\alpha h_\alpha )_\alpha$ converges to some $k \in G$ or tends to infinity. In the first case, we have
\[
  (\omega \bullet \Delta_\ast f)(g,h) = \lim_\alpha \Delta_\ast f(g g_\alpha, h_\alpha h) = \lim_\alpha  f(g g_\alpha h_\alpha h) = f(gkh) \qquad ((g,h) \in G \times H^\op)
\]
and in the second one
\[
   (\omega \bullet \Delta_\ast f)(g,h) = \lim_\alpha \Delta_\ast f(g g_\alpha, h_\alpha h) = \lim_\alpha  f(g g_\alpha h_\alpha h) = 0 \qquad ((g,h) \in G \times H^\op).
\]
In either case, $\omega \bullet \Delta_\ast f \in \SC(G \times H^\op)$ holds. This proves (iii).
\end{proof}
\section{Connes-amenability of $M(G)$}
Let $G$ be a locally compact group. As a consequence of Theorem \ref{ucthm}(iii), $\Delta_{M(G)}$ extends to an $M(G)$-bimodule homomorphism
$\Delta_{0,w^\ast} \!: \LUCSC(G \times G^\op)^\ast \to M(G)$:
\begin{proposition} \label{diagprop}
Let $G$ be a locally compact group such that $M(G)$ is Connes-amenable. Then there is ${\rm M} \in \LUCSC(G \times G^\op)^\ast$ such that
\[
  \mu \cdot {\rm M} = {\rm M} \cdot \mu  \quad (\mu \in M(G)) \qquad\text{and}\qquad \Delta_{0, w^\ast} {\rm M} = \delta_e.
\]
\end{proposition} 
\begin{proof}
Define a derivation
\[
  D \!: M(G) \to \LUCSC(G \times G^\op)^\ast, \quad \mu \mapsto \mu \tensor \delta_e - \delta_e \tensor \mu.
\]
It is easy to see that $D$ is $w^\ast$-continuous and attains its values in $\ker \Delta_{0,w^\ast}$. Being the kernel of a $w^\ast$-continuous bimodule homomorphism,
$\ker \Delta_{0,w^\ast}$ is a $w^\ast$-closed submodule of the normal, dual Banach $M(G)$-module $\LUCSC(G \times G^\op)^\ast$ and thus a normal, dual Banach $M(G)$-module
in its own right. Since $M(G)$ is Connes-amenable, there is thus ${\rm N} \in \ker \Delta_{0,w^\ast}$ such that
\[
  D \mu = \mu \cdot {\rm N} - {\rm N} \cdot \mu \qquad (\mu \in M(G)).
\]
The element
\[
  {\rm M} := \delta_e \tensor \delta_e - {\rm N}
\]
then has the desired properties.
\end{proof}
\begin{remark}
Since $\LUCSC(G \times G)^\ast$ is only a quotient of $\SC(G \times G)^\ast$, Proposition \ref{diagprop} does not allow us to conclude that 
$M(G)$ has a normal, virtual diagonal.
\end{remark}
\begin{lemma} \label{multipl}
Let $G$ and $H$ be locally compact groups. Then $\LUCSC(G \times H^\op)$ is an essential ideal of $\LUC(G \times H^\op)$.
\end{lemma}
\begin{proof}
Let $f \in \LUCSC(G \times H^\op)$, let $F \in \LUC(G \times H^\op) \subset \LUC(G \times H^\op)$, and let $\omega \in (G \times H^\op)_{\LUC}$. Let $( (g_\alpha, h_\alpha) )_\alpha$ be a net in
$G \times H^\op$ converging to $\omega$. Since
\[
  \omega \bullet (fF) = \lim_\alpha R_{(g_\alpha, h_\alpha)}(fF) = \lim_\alpha (R_{(g_\alpha, h_\alpha)}f)(R_{(g_\alpha, h_\alpha)}F) = (\omega \bullet f)(\omega \bullet F)
\]
with pointwise convergence on $G \times H$ and since $\omega \bullet F \in \LUC(G \times H^\op)$, it follows that $\omega \bullet (fF) \in \SC(G \times H)$. Hence, $\LUCSC(G \times H^\op)$ is an ideal of 
$\LUC(G \times H^\op)$. Since ${\cal C}_0(G \times H) \subset \LUCSC(G \times H^\op)$, it is even an essential ideal.
\end{proof}
\begin{theorem} \label{thm}
For a locally compact group $G$, the following are equivalent:
\begin{items}
\item $G$ is amenable.
\item $M(G)$ is Connes-amenable.
\end{items}
\end{theorem}
\begin{proof}
(i) $\Longrightarrow$ (ii): By \cite[Theorem 2.5]{Joh1}, $L^1(G)$ is amenable. Since $L^1(G)$ is $w^\ast$-dense in $M(G)$, \cite[Proposition 4.2]{Run} yields the Connes-amenability of $M(G)$.
\par
(ii) $\Longrightarrow$ (i): Let ${\rm M} \in \LUCSC(G \times G^\op)^\ast$ be as in Proposition \ref{diagprop}. View $\rm M$ as a measure on the
character space of the commutative $\cstar$-algebra $\LUCSC(G \times G^\op)$, so that $|{\rm M}| \in \LUCSC(G \times G^\op)^\ast$ can be defined in
terms of measure theory. It is routinely checked that $|{\rm M}| \neq 0$, and
\begin{equation} \label{yaeq}
  \delta_g \cdot | {\rm M} |  = | {\rm M} | \cdot \delta_g \qquad (g \in G). 
\end{equation}
By Lemma \ref{multipl} $\LUCSC(G \times G^\op)$ is an essential, closed ideal in $\LUC(G \times G^\op)$. We may therefore view $\LUC(G \times G^\op)$ as a $\cstar$-subalgebra of the multiplier algebra 
${\cal M}(\LUCSC(G \times G^\op))$. Since ${\cal M}(\LUCSC(G \times G^\op))$, in turn, embeds canonically into
$\LUCSC(G \times G^\op)^{\ast\ast}$, we may view ${\cal M}(\LUCSC(G \times G^\op))$ and thus $\LUC(G \times G^\op)$ as a $\cstar$-subalgebra of $\LUCSC(G \times G^\op)^{\ast\ast}$, 
so that, in particular, $\langle f, |{\rm M}| \rangle$ is well-defined for each $f \in \LUC(G \times G^\op)$. Note that the embedding of $\LUC(G \times G^\op)$ into
$\LUCSC(G \times G^\op)^{\ast\ast}$ is an $M(G)$-bimodule homormorphism (where the $M(G)$-bimodule action on $\LUC(G \times G^\op)$ is defined as on $\SC(G \times G)$).
Define
\[
  m \!: \LUC(G) \to \comps, \quad f \mapsto \langle f \tensor 1, | {\rm M}| \rangle.
\]
Since $f \tensor 1 \in \LUC(G \times G^\op)$ for each $f \in \LUC(G)$, $m$ is a well-defined, positive, linear functional. For $f \in \LUC(G)$ and $g \in G$, we have:
\begin{eqnarray*}
  \langle L_g f, m \rangle & = & \langle L_{(g,e)} (f \tensor 1), |{\rm M} | \rangle \\
  & = & \langle f \tensor 1, \delta_g \cdot |{\rm M}| \rangle \\
  & = & \langle f \tensor 1, |{\rm M}| \cdot \delta_g \rangle \\
  & = & \langle L_{(e,g)} (f \tensor 1), |{\rm M} | \rangle \\
  & = & \langle f \tensor 1, |{\rm M}| \rangle \\
  & = & \langle f, m \rangle.
\end{eqnarray*}
Normalizing $m$, we thus obtain a left invariant mean on $\LUC(G)$. Hence, $G$ is amenable by \cite[Theorem 1.1.9]{LoA}.
\end{proof}
\par
We believe that assertions (i) and (ii) in Theorem \ref{thm} are equivalent to:
\begin{items} \setcounter{enumi}{2} \it
\item $M(G)$ has a normal virtual diagonal.
\end{items}
Although we have been unable to prove this, Proposition \ref{compact} as well as the following corollary support this belief:
\begin{corollary} \label{disccor}
Let $G$ be a discrete group. Then the following are equivalent:
\begin{items}
\item $G$ is amenable.
\item $\ell^1(G)$ is Connes-amenable.
\item There is a normal, virtual diagonal for $\ell^1(G)$.
\end{items}
\end{corollary}      
\begin{proof}
(i) $\Longrightarrow$ (iii): If $G$ is amenable, $\ell^1(G)$ is amenable, so that there is a virtual diagonal ${\rm M} \in (\ell^1(G) \Tensor \ell^1(G))^{\ast\ast}$ for $\ell^1(G)$. Let $\rho \!: (\ell^1(G) \Tensor \ell^1(G))^{\ast\ast}
\to {\cal L}_{w^\ast}^2(\ell^1(G), \comps)^\ast$ denote the restriction map. Then $\rho({\rm M})$ is a normal, virtual diagonal for $\ell^1(G)$.
\par
Since (i) $\Longleftrightarrow$ (ii) by Theorem \ref{thm}, and since (iii) $\Longrightarrow$ (ii) for any dual Banach algebra, this proves the corollary.
\end{proof}                          
\begin{remark}
Since discrete groups are trivially inner amenable, the equivalence of (i) and (ii) in Corolllary \ref{disccor} can alternatively be deduced from \cite{LP}: If $\ell^1(G)$ is Connes-amenable, then so is $\VN(G)$ by
\cite[Proposition 4.2]{Run}, which, by \cite{LP}, establishes the amenability of $G$.
\end{remark}
\dated
\vfill
\begin{tabbing}
{\it Address\/}: \= Department of Mathematical and Statistical Sciences \\
\> University of Alberta \\
\> Edmonton, Alberta \\
\> Canada T6G 2G1 \\[\medskipamount]
{\it E-mail\/}: \> {\tt vrunde@ualberta.ca} \\[\medskipamount]
{\it URL\/}: \> {\tt http://www.math.ualberta.ca/$^\sim$runde/}   
\end{tabbing} 
\end{document}